\documentclass[12pt]{article}
\usepackage{amsmath,amssymb,amsfonts,amsthm}
\usepackage{eucal}%
\usepackage[all]{xy}
\usepackage{todonotes}
\newcommand {\qe} {\mathfrak{q}}

\newcommand {\al} {\alpha}

\newcommand {\lam}  {\lambda}

\newcommand{\bL}{\boldsymbol{\Lambda}}
\newcommand{\bal}{\boldsymbol{\alpha}}
\newcommand{\bphi}{\boldsymbol{\phi}}
\newcommand{\Kosz}{\textsf{Kosz}}
\newcommand{\ii}{\rm i}
\newcommand{\xmu}{\mathsf{x}_\mu}
\newcommand{\xmbox}{\mathsf{x}_\square}
\newcommand{\C}{\mathbb{C}}

\newcommand{\R}{\mathbb{R}}

\newcommand{\Z}{\mathbb{Z}}

\newcommand{\bR}{\mathsf{R}}
\newcommand{\cO}{\mathcal{O}}
\newcommand{\cZ}{\mathcal{Z}}

\newcommand{\cM}{\mathcal{M}}
\newcommand{\cE}{\mathcal{E}}

\newcommand{\bW}{\mathsf{W}}
\newcommand{\bV}{\mathsf{V}}

\newcommand{\cI}{\mathcal{I}}
\newcommand{\cJ}{\mathcal{J}}

\newcommand{\cN}{\mathcal{N}}
\newcommand{\cX}{\mathcal{X}}
\newcommand{\bE}{\mathsf{E}}

\newcommand{\Pp}{\mathbf{P}^1}
\newcommand{\Pt}{\mathbf{P}^2}

\newcommand{\zz}{\mathfrak{z}}
\newcommand{\ww}{\mathsf{g}}

\newcommand{\bT}{\mathsf{T}}
\newcommand{\cF}{\mathcal{F}}
\newcommand{\cG}{\mathcal{G}}
\newcommand{\cH}{\mathcal{H}}
\newcommand{\cS}{\mathcal{S}}
\newcommand{\fZ}{\mathfrak{Z}}

\newcommand{\lang}{\left\langle}
\newcommand{\rang}{\right\rangle}

\DeclareMathOperator{\Hilb}{Hilb}
\DeclareMathOperator{\Ext}{Ext}
\DeclareMathOperator{\Hom}{Hom}

\DeclareMathOperator{\tr}{tr}

\DeclareMathOperator{\pt}{pt}

\newtheorem{Theorem}{Theorem}
\newtheorem{Lemma}{Lemma}
\newtheorem{Proposition}[Lemma]{Proposition}
\newtheorem{Corollary}[Lemma]{Corollary}
 
\begin{document}

\title{Five dimensional gauge theories\\ and\\ Vertex operators}
\author{Erik Carlsson$^1$, Nikita Nekrasov$^{1,2,3,4}$, and Andrei Okounkov$^{3,5}$}

\date{June 26, 2013} 
\maketitle
\begin{center}
${}^{1}$ Simons Center for Geometry and Physics,
Stony Brook NY 11794-3636 USA\\

${}^{2}$  Institut des Hautes Etudes Scientifiques,  Bures-sur-Yvette 91440 France\\

$^{3}$ Kharkevich Institute for Information Transmission Problems,
Lab. 5, Moscow 127994 Russia
\footnote[0]{On leave of absense}\\

$^{4}$ Alikhanov Institute of Theoretical and Experimental Physics,
Moscow 117218 Russia\footnotemark[0]\\
$^{5}$ Department of Mathematics, Columbia University, New York USA
\end{center}
\bigskip
\begin{center}
\it To Grigori Olshanski, on his 64th birthday 
\end{center}

\begin{abstract}
We study supersymmetric gauge theories in five dimensions, using their relation to the K-theory of the moduli spaces of torsion free sheaves. In the spirit of the BPS/CFT correspondence the partition function and the expectation values of the chiral, BPS protected observables are given by the matrix elements and more generally by the correlation functions in some $q$-deformed conformal field theory in two dimensions. We show that the coupling of the gauge theory to the bi-fundamental matter hypermultiplet inserts a particular vertex operator in this theory. In this way we get a generalization of the main result of \cite{CO} to $K$-theory. The theory of interpolating Macdonald polynomials is an important tool in our construction. 
\end{abstract}

\section{Introduction}

\subsection{}

We begin with a brief explanation of the historical motivation behind the problems 
studied in this paper. The reader may want to consult, for example, \cite{N,NO,Oicm} for more 
details.

In \cite{N}, partition functions of certain $4$-dimensional supersymmetric quantum gauge theories
were given a mathematical definition as equivariant integrals over the instanton moduli spaces. 
The integrands depend on the matter content of the theory and 
represent characteristic classes of natural vector bundles over instanton moduli. 

This description can be generalized \cite{Nekrasov:1996cz} to $5$-dimensional gauge theories on a product of the 
original $4$-fold $X$ and a circle ${\cS}^1$. Integration in equivariant cohomology is then 
replaced by push-forwards in equivariant $K$-theory. 

Since the early days of the studies of $S$-duality in gauge theory \cite{Vafa:1994tf} and later in string theory the  rich structure of conformal field theories in two dimensions was suspected to be within the realm of the four dimensional supersymmetric gauge theories. Partly this feeling was fueled by the constructions in \cite{Nak1,Nak2,Nak} which realized affine Kac-Moody algebras and their deformations acting via some correspondences on the cohomology and the K-theory of the moduli spaces of gauge instantons on ${\bf R}^{4}$ and resolutions of its orbifolds. 
 
 It was proposed in \cite{N} to interpret the chiral fields of the two dimensional conformal field theory describing the BPS sector of the four dimensional ${\mathcal N}=2$ supersymmetric theory, as well as its equivariant generalizations, the so-called $\Omega$-deformation, as the zero modes of the chiral tensor field in the theory on the M-theory or NS5 brane, which engineers the  ${\mathcal N}=2$ theory, by wrapping on a Riemann surface, the Seiberg-Witten curve $\Sigma$. It was later realized in \cite{Losev:2003py, NO} that the conformal field theory need not to live on the Seiberg-Witten curve, but rather some auxiliary curve $C$. The Seiberg-Witten curve $\Sigma$, an emergent geometry, is a ramified covering of $C$. 
 This circle of ideas led to the concept of the  BPS/CFT correspondence \cite{NA}. In recent years several interesting examples of the BPS/CFT correspondence were found. In \cite{Witten:1997sc, Gaiotto:2009we, Gaiotto:2009hg} a class of four dimensional theories, the so-called $S$-class theories, labelled by a Riemann surface $C$ together with a choice of its pants decomposition, and a choice of the $ADE$ simple Lie group, was proposed. In \cite{Alday:2009aq} the evidence for the equality of the partition functions of the $\Omega$-deformed $S$-class theories of $A_1$  type with the conformal blocks of Liouville theory on $C$ was given, leading to the AGT conjecture. The general $A_k$ case is conjectured to be related to the Toda conformal field theories.  
 In \cite{NS} a particular limit of the general $\Omega$-deformed theory was conjectured to capture the spectrum of a quantum integrable system, which for the $S$-class four dimensional gauge theories is a version of quantum Hitchin system. The investigation of the five dimensional theories compactified on a circle should lead to a relativistic \cite{Nekrasov:1996cz} (difference) version of the quantum Hitchin system, which is yet another reason to be interested in the story below.

\subsection{}

In this paper we shall focus on the theories with the gauge groups
which are the products of the unitary groups, say $U(r) \times U(r')$,
with the matter in the bifundamental representation. We study the
theory on $X \times {\cS}^1$ where $X$ is a complex surface, for
example  ${\C}^{2}$. 
In holomorphic description, the
$U(r)$ instantons correspond to rank $r$ holomorphic bundles on $X$. More precisely, 
the moduli space of instantons is partly compactified to the moduli space $\cM(r)$ of the 
torsion free rank $r$ sheaves $\cE$ on a certain compactification $\bar X$ of $X$, with a trivialization 
\begin{equation}
{\bf\Phi} : {\cE} \vert_{D} \longrightarrow^{\kern -.2in \approx\, } \ {\cO}_{D}^{{\oplus} r}
\label{eq:trivphi}
\end{equation}
over the compactification divisor $D = {\bar X} \backslash X$. 

On
the product $\cM(r) \times \cM(r')$ of instanton moduli, there is a
natural sheaf formed by  $\Ext$-groups between the bundles in questions
(see Section \ref{s_Ext} for precise definitions). It describes 
matter in the bifundamental representation. 

As any $K$-theory class on the product, the extension bundle defines
a Fourier-Mukai operator
\begin{equation}
\Phi_{\Ext} : K_{G \times T} (\cM(r)) \to K_{G\times T} (\cM(r')) 
\label{eq:fmop}
\end{equation}
where $G$ is the group of constant gauge transformations which act on
$\cM (r)$ by changing the trivialization, ${\bf\Phi} \mapsto g \circ
{\bf\Phi}$, and $T = {\C}^{\times} \times {\C}^{\times}$ is the
maximal torus of the group $GL(2, {\C})$. It acts by  complex rotations of $X$.

\subsection{}

The structure of $G \times T'$-equivariant cohomology and $G \times T'$-equivariant $K$-theory of $\cM (r)$
greatly simplifies when $T'$ preserves a holomorphic symplectic form on $X$. 
In this case, it was noted in  \cite{NO} that $\Phi_{\Ext}$ is a vertex
operator for Nakajima's Heisenberg algebra, see also  \cite{LS, Losev:2003py}.

The correct fully $G \times T$-equivariant generalization was found for $r=r'=1$ in 
\cite{ErikTh,CO}. Here we generalize this result to $K$-theory. 

The definition of the operator $\Phi_{\Ext}$ is given in Section \ref{s_Phi}
and our formula for it in terms of certain deformed Heisenberg 
operators is the Corollary \ref{C1} in Section \ref{s_C1}. 

\subsection{}
 
Equivariant $K$-theory does not permit a number of dimension-based arguments used in \cite{CO}. 
Therefore, in the present paper we take a very different approach. 
The computation of $\Phi_{\Ext}$ is eventually reduced to  a certain fundamental 
fact in the harmonic analysis for Cherednik algebra, namely the so-called 
Macdonald-Mehta-Cherednik identity \cite{C2,EK}. 

This identity computes the values of the Hermitian form
\begin{equation}
(f,g) \mapsto \lang f, \Theta \, g\rang_\Delta
\end{equation}
in the basis of Macdonald polynomials. Here $\lang \, \cdot \,, \, \cdot \,\rang_\Delta$ 
is Macdonald Hermitian inner product associated to a root system and $\Theta$ is 
the theta function of the corresponding 
weight lattice. Our formula for $\Phi_{\Ext}$ is 
obtained from the $GL(N)$ case as $N\to \infty$. In this limit, the 
vertex operators arise as a stabilization of the triple-product factorization of $\Theta$. 

The same problem can be also approached from the angle of interpolation Macdonald
polynomials, see Section \ref{s_interp}. 
 For brevity, in this note we use it in reverse and deduce an interesting 
formula for interpolation polynomials from the Macdonald-Mehta-Cherednik identity, see 
Proposition \ref{NGa}. 

\subsection{}

This paper was written in 2008, in time for Grigori Olshanski's 
60th birthday, but had only a very limited
circulation. It was meant to be a part 
of a larger project which, regrettably, has been slow to materialize. 
We thought that this could be a good occasion to revisit our old
manuscript. 

While there has been a great progress since 2008 in the understanding 
of the issues considered here (see for example \cite{Aetal} and also
\cite{IKS} 
where an equivalent formula may be found), we feel that 
our basically elementary reduction to Macdonald-Mehta-Cherednik formulas,
as well as the connection to interpolation polynomials,
 could be of interest. 

\subsection{}

The interpolation Schur functions and, more generally, 
Macdonald polynomials, a subject pioneered by G.~Olshanski, 
have been finding numerous applications in representation theory, 
algebraic geometry, and mathematical physics, the present paper 
being yet another example. 

An important question, on which we touch here, is why a beautiful
theory of such polynomials exists for $GL(n)$ and also $BC_n$, 
but not, for example, for the roots system 
$A_n$. In section \ref{s_interp}, we 
discuss a very basic reason for this: the weight lattice $\Z^n$ 
of $GL(n)$ is an orthogonal direct sum of $\Z$'s. 

\subsection{} 

NN was partially supported by RFBR grant 12-02-00594, and by the
Agence Nationale de la Recherche under grant 12-BS05-003-01.  
AO thanks NSF for finantial support under FRG grant 0853560.

\section{Background}

An important player in what follows is the moduli space $\cM (1) = \amalg_{n} \Hilb_n$ which 
has different incarnations. For fixed $n \geq 0$ it is a moduli space of $n$ D0 branes bound to a single D4 brane in the IIA string theory \cite{Douglas:1995bn}, it is a moduli space of $U(1)$ instantons of charge $n$ on a noncommutative Euclidean space ${\R}^4$ \cite{Nekrasov:1998ss}, it is a resolution of singularities of the moduli space of charge $n$ $U(r)$ instantons on ${\R}^{4}$ for $r=1$, \cite{Nak2}, and it is a Hilbert scheme of $n$ points on $\C^2$. The larger family 
of varieties that includes moduli $\cM(r)$ of torsion-free
$\C[x,y]$-modules of rank $r>1$ and more general Nakajima
varieties also play an important role in supersymmetric 
gauge theories and the question that we study here for 
Hilbert schemes is interesting and important in 
that greater generality. In this paper, however, we 
focus on Hilbert schemes. Thus, from the gauge theory 
perspective, we restrict ourselves to abelian gauge 
theories on $\C^2\cong \mathbb{R}^4$.

\subsection{Hilbert schemes of points} 

\subsubsection{}

For $n=0,1,2,\dots$, let $\Hilb_n$ denote the Hilbert scheme of $n$ points in the 
plane $\C^2$. By definition, a point $I \in \Hilb_n$ is 
an ideal $I\subset \C[x,y]=\cO_{\C^2}$ such that 
$$
\dim_\C \C[x,y] /I = n \,.  
$$  
It has a natural structure of complex algebraic variety 
and is known to be smooth and irreducible of dimension 
$2n$, see for example \cite{Nak} for an introduction. 

\subsubsection{}

There is a universal subscheme 
$$
\fZ \subset \Hilb_n \times \C^2
$$
such that its intersection with $\{I\}\times \C^2$, 
$I\in\Hilb_n$,  is the subscheme of $\C^2$ defined by 
the ideal $I$. The sheaf of ideals $\cI$ of $\fZ$ is called
the universal ideal sheaf on $\Hilb_n \times \C^2$.

\subsubsection{}

The group $GL(2)$ acts naturally on $\C^2$ and, hence, on $\Hilb_n$. 
In this paper, we will be working 
with the $GL(2)$-equivariant $K$-theory of $\Hilb_n$. We also pick 
a maximal torus $T$ in $GL(2)$. 

We denote by $q^{-1},t^{-1}$ the torus weights on $\C^2$, 
that is 
$$
[\C^2] = q^{-1} + t^{-1} 
$$
in the representation ring of $G$. This means that torus weights in $\C[x,y]$
are $\{q^i t^j\}$, $i,j\ge 0$. 

\subsubsection{} 

The physics background for this problem is the following. Consider ${\mathcal N}=1$ supersymmetric gauge theory in five dimensions, and compactify it on a circle $S = {\cS}^{1}$, of circumference $2\pi r$, 
with twisted boundary conditions, 
such that the remaining space ${\R}^{4}$ is rotated by an element $(q,t)$ of the maximal torus of the spin cover $Spin(4)$ of the rotation group. In addition the fermions
in the vector multiplet are subject to the additional $SU(2)$ $R$-symmetry twist 
\begin{equation}
\left( 
\begin{matrix}
(qt)^{\frac 12} & 0 \\
0 & (qt)^{-\frac 12} \end{matrix} 
\right) 
\label{eq:rsym}
\end{equation}
which is chosen so as to preserve some fraction of the original supersymmetry. 
The preserved supercharge can be interpreted \cite{Nekrasov:1996cz, Baulieu:1997nj} as the equivariant de Rham differential
acting in the space of equivariant forms on the loop space on the space of four dimensional gauge fields. The source of the loops is the circle $S$ of compactification of the five dimensional theory down to four dimensions. 

Actually, the construction described above would lead to $|q|=|t|=1$. However, by turning on some additional fields in the background of $\cN =2$ supergravity one can make $q,t$ general complex numbers. One can also describe this background (the general $\Omega$-background) by lifting the theory to six dimensions, and then compactifying on a two-torus $T^2$. The flat $Spin(4)$ connection on $T^2$, in the limit where the torus collapses to a circle, becomes, effectively, a flat $Spin(4, {\C})$-connection, which can always be reduced  to a $T$-connection. The pair $(q,t)$ is its holonomy.

The common lore states that the gauge theory in five dimensions is not well-defined as the quantum field theory unless it has a cutoff. However, there exist nontrivial fixed points at the strong coupling \cite{Seiberg:1996bd}, which can be realized as limits of compactifications
of M-theory on Calabi-Yau three folds, or via brane configurations. Another possible ultraviolet completion of the theory is given by its embedding in the $(2,0)$ superconformal theory compactified on a circle $S^\prime$. 

The $(2,0)$ theory is, in turn, a limit of the theory on a single M5 brane. The
worldvolume $\cX$ of the M5 brane is six dimensional, and the normal bundle 
to $\cX$ is a rank $5$ vector bundle $\cN$ over $\cX$. 
The rotation \eqref{eq:rsym} acts geometrically  by rotating
the fibers of $\cN$. A generic rotation of a five dimensional Euclidean space 
has an invariant
one dimensional line $\R$. 

The full setup therefore involves M-theory on the eleven-dimensional manifold, 
which is a product of a real line $\R$, a circle $S^\prime$
 and a rank $8$ vector bundle $N \times X$ over $S$. 
 The geometry
 of the vector bundle is determined by the two rotation angles which translate to $q$ and $t$ parameters of our story.

\subsubsection{}

While $\Hilb_n$ is not compact, the Euler characteristic (a virtual representation of $GL(2)$)
\begin{equation}
\chi_{\Hilb}(\cF,\cG) = \sum_{i} (-1)^{i} \Ext_{\Hilb}^i(\cF,\cG)\label{chi}
\end{equation}
is well-defined as a $GL(2)$-module for any pair $\cF,\cG$ of 
coherent sheaves on $\Hilb_n$. Namely, it is 
a direct sum of irreducibles with finite multiplicities. 

To see this, one can use
the natural map 
$$
\pi: \Hilb_n \to (\C^2)^n / S(n) 
$$
that takes an ideal $I$ to the support of $\C[x,y]/I$, 
counting multiplicity. The map $\pi$ is proper. 
On the other hand, the center of $G$ acts on $\C^{2n}/S(n)$
with positive weights.

\subsection{Matrix integrals}

Using the ADHM construction of the moduli spaces of instantons and torsion free sheaves, the Euler characteristics  \eqref{chi} can be written as an integral over $n\times n$ matrices (see \cite{Lossev:1997bz, Moore:1997dj, Moore:1998et} for similar integrals). This integral can be further reduced, using supersymmetric localization, to the $n$-fold contour integral:
\begin{equation}
\begin{aligned}
& {\chi}_{\Hilb_n} (\cF,\cG) = \\
& \qquad \frac{1}{n! \, (2 \pi {\ii})^n} \oint \prod_{i=1}^{n} 
\frac{dx_i / x_i}{T(w, x_i; qt)} \, f(x) g(x^{-1})\,  {\tilde\Delta}_{q, t} (x)\, \\
& \qquad\qquad\qquad {\tilde\Delta}_{q,t} (x) =   {\Gamma}^{\prime}(1)^{n} \, \prod_{1 \leq i \neq j \leq n} {\Gamma}(x_i /x_j ) \, , \\
& \qquad\qquad\qquad\qquad {\Gamma}(y) = \frac{(1-y)(1-qt y)}{(1-q y)(1-t y)}\, , \\
& \\
& \qquad\qquad\qquad\qquad\qquad T(w, y; Q^2) = (1 - w Q/y)(1 - y Q/w) 
\end{aligned}
\label{eq:chifg}
\end{equation}
where the contours are the circles $| x_{i} | = 1$,  $|w| \sim 1$ and we assumed $|q|, |t | \ll 1$. 
Finally, the $K$-theory classes $\cF$ and $\cG$ are represented by the symmetric functions $f(x)$ and $g(x)$, with the dual bundle $\cG^{\vee}$ corresponding to the symmetric function $g^{\vee}(x) = g(x^{-1})$. The K-theory of $\Hilb_n$ is generated by the tensor functors of the tautological rank $n$ vector bundle $V = {\cO}/I$. Write the Chern character of $V$ as
\begin{equation}
Ch(V) = \sum_{i=1}^{n} x_{i}
\end{equation} 
An irreducible representation $\mu$
of $GL(n)$ corresponding to a length ${\ell}({\mu}) \leq n$ partition $\mu$ gives rise to the
associated vector bundle $V_{\mu}$ over $\Hilb_n$, which can be constructed using, e.g. the Young symmetrizer. Its Chern character $Ch(V_{\mu})$ is represented by the Schur function
$s_{\mu}(x)$:
\begin{equation}
s_{\mu}(x) = \frac{{\rm Det}_{i,j} \Vert x_{i}^{{\mu}_{j} - j + n} \Vert}{ {\rm Det}_{i,j} \Vert x_{i}^{- j + n} \Vert}  
\end{equation}
The contour integral \eqref{eq:chifg} has to be supplemented with the prescription, that it can be computed by residues, the poles corresponding to the $T$-fixed points in $\Hilb_n$. The latter are labeled by the size $|{\lam}| = n$ partitions $\lam$. The possible poles coming from the singularities of $g(x^{-1})$ for $|x_{i}|< 1$
should be dropped. 

The pole corresponding to the partition $\lam$ is at (up to a $S(n)$-permutation):
\begin{equation}
x_{\square} = w \, q^{i_\square - \frac 12} t^{j_\square - \frac 12}, \qquad \square \in \lam
\label{eq:cont}
\end{equation}
where $\square = (i_\square, j_\square)$ with $1 \leq i_\square \leq {\ell}({\lam})$, 
$1 \leq j_{\square} \leq {\lam}_{i_{\square}}$.

\subsection{From matrix to Macdonald}

Note that a symmetric function of $x_i$'s can be expressed as a function of the power-sums
$p_k = \sum_{i=1}^{n} x_i^k$. One can take, for example, $p_k$ with $1 \leq k \leq n$, as
the generators. 
The measure ${\tilde\Delta}_{q,t}$ in \eqref{eq:chifg} can be reexpanded:
\begin{equation}
{\tilde\Delta}_{q,t}(x) = {\exp}\,  \sum_{k=1}^{\infty} \frac 1k \left( - (1-q^{k})(1-t^{k})p_{k}p_{-k}  + n \right)
\label{eq:macdo1}
\end{equation}
We can also expand the $T$-factors in \eqref{eq:chifg}
\begin{equation}
\prod_{i=1}^{n} T(w, x_i ; Q^2 )^{-1} = {\exp} \, \sum_{k=1}^{\infty} \frac{Q^k}{k} \left( 
p_{k}  w^{-k} + p_{-k} w^{k}  \right) 
\end{equation}
Now define, for $k > 0$:
\begin{equation}
{\al}_{k} = \frac{(w\sqrt{qt})^{k}}{1-q^{k}} - ( 1 - t^{k} )  p_{k}
\label{eq:alk}
\end{equation}
Evaluation of ${\al}_{k}$ on $(x_i)_{i=1}^{n}$ corresponding to a particular pole \eqref{eq:cont}
gives:
\begin{equation}
{\al}_{k} = (w\sqrt{qt})^{k} \sum_{i=1}^{\infty} q^{k(i-1)} t^{k{\lam}_{i}}  = 
\sum_{i=1}^{\infty} {\xi}_{i}^{k}
\label{eq:alkplam}
\end{equation}
Now the idea is to pass from the variables $x_i$, which have the meaning of the
positions of D0 branes on $S \times {\R}$, to the variables $\xi_i$. 
The formula \eqref{eq:alk} (but not \eqref{eq:alkplam}) makes sense for all
$k \neq 0$. In terms of $\al_k$ the product of the measure factors in \eqref{eq:chifg}
reads simply as:
\begin{equation}
(qt)_{q , t} {\tilde\Delta}_{q,t}(x)  \prod_{i=1}^{n} T(w, x_i ; qt )^{-1}
= {\exp}\,  \sum_{k=1}^{\infty}
\frac 1k \left( n  - \frac{1-q^{k}}{1-t^{-k}}\,  {\al}_{k}{\al}_{-k}  \right)
\label{eq:chixi}\end{equation}
where 
\begin{equation}
(x)_{q_{1},q_{2}}  = \prod_{i,j=0}^{\infty} ( 1 - x q_{1}^{i}q_{2}^{j} ) 
\label{eq:xqt}
\end{equation}
 Now, if instead of \eqref{eq:alkplam} with an infinite number of ${\xi}_{i}$'s 
 we take
 a finite number $N$ of $\xi_{i}$'s, then, up to a divergent constant, \eqref{eq:chixi} defines a measure
 on symmetric functions $F({\xi})$, and a hermitian form on symmetric polynomials
 $P({\xi})$:
 \begin{equation}
 \begin{aligned}
 & \langle P, Q \rangle = \frac{1}{(2\pi {\ii})^N N!} \oint \prod_{a=1}^{N} \frac{d{\xi}_{a}}{{\xi}_{a}} \, P({\xi}) Q({\xi}^{-1})  \, {\Delta}_{q,t}({\xi}) \\
&  \qquad {\Delta}_{q,t}({\xi}) = {\exp}\,  \sum_{k=1}^{\infty}
\frac 1k  \frac{1-q^{k}}{1-t^{-k}}\,  \left( N  - {\al}_{k}{\al}_{-k}  \right)
= \\
& \qquad\qquad\qquad\qquad =  \prod_{n=1}^{\infty} \prod_{1\leq a \neq b \leq N}
\frac{1 - t^{n} {\xi}_{a}/{\xi}_{b}}{1 - q t^{n} {\xi}_{a}/{\xi}_{b}} \end{aligned}
\label{eq:macdom}
\end{equation} 
 This is the scalar product leading to Macdonald polynomials \cite{Mac}.

\subsection{The Ext bundle} \label{s_Ext}

\subsubsection{}

The tangent bundle to $\Hilb_n$ may be expressed in terms of the universal 
sheaf $\cI$ as follows:
\begin{align}
T\Hilb &= \chi_{\C^2}(\cO) - \chi_{\C^2}(\cI,\cI) \,, \label{THilb}\\
  &= \chi_{\C^2} (\cO_\fZ,\cO) + \chi_{\C^2} (\cO,\cO_\fZ) - 
\chi_{\C^2} (\cO_\fZ,\cO_\fZ) \notag 
\end{align}
where the subscripts mean that we take the 
push-forward along the $\C^2$-factor,
$$
\cO_\fZ = \cO/\cI
$$
is the structure sheaf of the universal subscheme, 
and the push-forward in the second line of 
\eqref{THilb} is well defined because $\fZ$ is proper
in the $\C^2$-direction.

\subsubsection{}

It is natural to generalize \eqref{THilb} to a bundle 
$\bE$ on the product of two Hilbert schemes with fiber 
\begin{align}
  \label{defE}
 \bE(I,J) &= \chi_{\C^2}(\cO) - \chi_{\C^2}(I,J) \\
    &=  \Ext^1_{\Pp}(I,J(-1)) \notag 
\end{align}
over a point 
$$
(I,J)\in \Hilb_n \times \Hilb_m \,. 
$$
Here $(-1)$ denotes the twist by the line $\Pt\setminus \C^2$ at
infinity of $\C^2$.

In gauge theory, a bundle like $\bE$ is related to matter in 
bifundamental representation of two gauge groups. 
Here we have the $U(1)\times U(1)$ bifundamental matter. 

\subsubsection{}
Any $K$-theory class $\cE$ on the product of two 
varieties $X \times Y$, defines an operator 
$$
\Phi_\cE : K(X) \to K(Y)
$$
by the formula
$$
\Phi_\cE(\cF) = {p_X}_*  \left( \cE \otimes p_Y^*\cF \right)\,.
$$
Here 
$$
p_X : X \times Y \to Y
$$
is the projection \emph{along} $X$, and similarly for $p_Y$. 

\subsubsection{}\label{s_Phi}
Our goal in this paper is to describe the operators 
$$
\Phi_{\Lambda^k\bE}\,, \quad k=1,2, \dots\,,
$$
where $\Lambda^k\bE$ denote the exterior powers of $\bE$. It will be 
convenient to package them into a generating function 
\begin{equation}
  \label{defPhi(m)}
  \bW(m) = \sum (-m)^k \Phi_{\Lambda^k \bE} \,. 
\end{equation}
The argument $m$ of this generating function is related to the mass $\mu$
of the bifundamental matter in gauge theory, via
\begin{equation}
m = e^{r{\mu}}
\label{eq:mass}
\end{equation}

\subsubsection{}

Concretely, we want to describe them as operators on symmetric
functions using the identification 
$$
K_{GL(2)} (\Hilb_n) \cong \textup{(symmetric functions of degree $n$)} \otimes K_{GL(2)}(\pt)\,,
$$
that follows from the work of Haiman \cite{Hai1} and Bridgeland, King, and Reid \cite{BKR},
see below. 
 
In  terms of symmetric functions, this task may be phrased as computing matrix
coefficients of a certain operator in the basis of Macdonald polynomials.

\subsection{Hermitian inner product in $K$-theory}

\subsubsection{}

For any $Z$, the Euler form \eqref{chi} on $K(Z)$ is sesquilinear, 
but is not symmetric or antisymmetric. This 
may be fixed if the canonical line bundle $K_Z$ has a square root $K_Z^{1/2}$. 
For example, for the Hilbert scheme, the canonical class is a pure character
$$
K_{\Hilb_n} = \left(\Lambda^2\C^2\right)^{-n} 
$$
and its square root may be extracted on a double cover of $GL(2)$. 

\subsubsection{}

One defines
\begin{equation}
\lang \cF,\cG \rang_Z = \chi_Z\left(\cF,K_Z^{1/2} \otimes \cG\right)\,, 
\quad \cF,\cG \in K(Z)  
\label{def_form}
\end{equation}
where $K_Z$ is the canonical line bundle of $Z$. By Serre duality 
$$
\lang \cF,\cG \rang_Z = (-1)^{\dim_\C Z} \lang \cG,\cF \rang_Z^\vee \,. 
$$
In this paper, we will meet only even-dimensional varieties, so 
the form \eqref{def_form} will be an Hermitian form for us. 

\subsubsection{}

We will take adjoints for maps between $K$-theories with respect to the 
form \eqref{def_form}. Because the form is symmetric, we don't need
to distinguish between left and right adjoints. In particular, one has
$$
\left(\Phi_\cE\right)^* = \Phi_{\cE_\textup{adjoint}}
$$
where
$$
\cE_\textup{adjoint} = \cE^\vee \otimes K_X^{1/2} \otimes K_Y^{1/2} \,. 
$$
Here $\cE^\vee$ is the dual bundle and the square roots $K_{X,Y}^{1/2}$ are pulled back to $X \times Y$ 
along the two projections.

\section{Factorization of the Ext operator}

\subsection{}

The first step in our computation is completely abstract and general. 
We consider
$$
\cI^{\boxtimes k} \in K_{GL(2)\times S(k)}\left(\Hilb \times (\C^2)^k\right). 
$$
This is $k$-th tensor power of $\cI$, except we take the 
$\otimes$-product along $\Hilb$ and the exterior $\boxtimes$-product
along $\C^2$. 

In other words, we consider $\Hilb \times \C^2$ as a variety over $\Hilb$
and take exterior tensor product in that category. 
The symmetric group $S(k)$ acts on $\cI^{\boxtimes k}$ permuting 
the factors.

\subsection{}

For any topological space $Z$ we define:
\begin{equation}
e^{Z} = {\rm pt} \amalg \, \bigcup_{k=1}^{\infty} Z^{k}/S(k)
\end{equation}

For any $Z$, we define $K(Z^k/S(k))$ as $K_{S(k)}(Z^k)$ with 
Euler form 
\begin{equation}
\chi_{Z^k/S(k)} = \left(\chi_{Z^k}\right)^{S(k)} 
\label{chi_orb}
\end{equation}
and set
$$
 K\left(e^Z\right)
 \overset{\textup{\tiny def}}= \bigoplus_{k\ge 0}  K(Z^k/S(k)) \,. 
$$

For any $\cF,\cG\in K(Z)$ we then have 
\begin{equation*}
\chi_{Z^k/S(k)}(\cF^{\boxtimes k},\cG^{\boxtimes k})  = 
S^k \chi_{Z}(\cF,\cG) \,. 
\end{equation*}
Here $Z$ can be a variety over any base so, in particular, 
we have 
\begin{equation}
\chi_{\C^{2k}/S(k)}(\cI^{\boxtimes k},\cJ^{\boxtimes k}) = 
S^k \chi_{\C^2}(\cI,\cJ) \,,\label{Sk}
\end{equation}
in $K_{GL(2)}(\Hilb)$.

\subsection{}

Define 
\begin{equation}
\bV: K_{GL(2)}(\Hilb) \to K_{GL(2)} \big(e^{\C^2}\big)
\label{defbV}
\end{equation}
by 
$$
\bV = \bigoplus_{k\ge 0} \Phi_{\cI^{\boxtimes k}} \,. 
$$
We also define the grading operators $L_0$ in 
in the domain and target of \eqref{defbV} as follows: it acts in 
$$
K_{GL(2)}(\Hilb_n)\,, \quad K_{GL(2)}(\C^{2n}/S(n))
$$
with the eigenvalue $n$. 

{}From \eqref{Sk} and definitions, we have 
$$
\bV^*  \left(\frac{m}{\sqrt{qt}}\right)^{L_0} 
\bV = \sum_k m^k \Phi_{S^k \chi_{\C^2}(\cI,\cJ)}\,,
$$
where the star $*$ denotes the adjoint with respect to \eqref{def_form}.

\subsection{} 

{}From
$$
\chi_{\C^2}(\cI,\cJ) = \chi_{\C^2}(\cO) - \bE(\cI,\cJ)
$$
we conclude
\begin{equation}
S^k \chi_{\C^2}(\cI,\cJ) = 
\sum_{l=0}^k (-1)^l \, S^{k-l} 
\chi_{\C^2}(\cO) \otimes \Lambda^l \, \bE \,. 
\label{S2wedge}
\end{equation}
Now the convenience of having a generating function \eqref{defPhi(m)} 
becomes clear. We compute, cf. \eqref{eq:xqt}
$$
\sum_{k\ge 0} m^k S^k \chi_{\C^2}(\cO) = 1/ (m)_{q,t} \,, 
$$
Therefore, formula \eqref{S2wedge} gives 
$$
\bW(m) = (m)_{q,t} \, 
\sum_k m^k \, \Phi_{S^k \chi_{\C^2}(\cI,\cJ)}  \,. 
$$
Thus we have proven the following 

\begin{Theorem}\label{T1}
We have 
$$
\bW(m) = (m)_{q,t} \, \bV^* \, 
\left(\frac{m}{\sqrt{qt}}\right)^{L_0}  \, \bV \,. 
$$
\end{Theorem}

Note that the argument was completely abstract and applies to
many more general moduli of sheaves. 

\section{Universal sheaf as an operator}

\subsection{}

Our next step is the identification of the operator $\bV$ 
defined in \eqref{defbV}. The description we seek goes via
the identification 
\begin{equation}
 \xymatrix{
K_{GL(2)}(\Hilb) \ar@{->}[r]^{\bV} 
 \ar[d]_{{\mathcal{B}} \sim} & K_{GL(2)}(e^{\C^2}) \ar[d]_{\sim} \\
\bL \ar@{->}[r]^{\bV} & \bL \\
} \label{comm1}
\end{equation}
where $\bL$ denotes symmetric functions 
over 
$$
\bR = \Z [ q,t] \left[\frac{1}{1-q^k t^l}\right] \,,
$$
and we similarly extend the scalars in the top row of \eqref{comm1}. 
The vertical 
identifications in \eqref{comm1} 
are as follows. 

\subsection{}

The map 
\begin{equation}
K_{{GL(2)} \times S(n)}(\C^{2n}) \to K_{{GL(2)}\times S(n)}(\pt) \cong \bL_n
\label{global_sec}
\end{equation}
is given simply by taking a global section. Here $\bL_n$ denotes the space of 
symmetric functions of degree $n$ and the isomorphism with 
$K_{{GL(2)}\times S(n)}(\pt)$ takes a Schur function $s_\lambda$ to 
the corresponding 
irreducible representation $[\lambda]$ of the symmetric group 
(viewed  as a trivial ${GL(2)}$-module). In other words, 
$$
s_\lambda \leftrightarrow [\lambda] \otimes \cO_0
$$
where $\cO_0$ is the skyscraper sheaf at the origin $0\in \C^{2n}$.

\subsection{} 
The map \eqref{global_sec} is an isomorphism of algebras
with respect to induction of representations. Concretely, 
the product of two equivariant sheaves 
$$
\cF_i \in K_{{GL(2)}\times S(n_i)}(\C^{2n_i}) 
$$
is defined to be 
$$
\cF_1 \cdot \cF_2 = \C S(n_1+n_2) \otimes_{\C S(n_1) \otimes \C S(n_2)} 
\cF_1 \boxtimes \cF_2 \,. 
$$

\subsection{} 
The map \eqref{global_sec} does not preserve $\Hom$'s, that is, 
it induces a new Hermitian form on $\bL$. It is best described
in terms of the tensor inverse 
$$
\textsf{Kosz} = \bigoplus (-1)^i \Lambda^i \left( \C^{2n} \right)^\vee \,,
$$
of the polynomial representation $S^*\!\left( \C^{2n} \right)^\vee$. 
It appears in the Koszul resolution 
$$
\Kosz \otimes \cO_{\C^{2n}}  \to \cO_0
$$
of the skyscraper sheaf at the origin. The Hermitian form 
on induced from the orbifold $\exp(\C^2)$ is 
$$
\lang M, M' \rang = (qt)^{n/2} \, \Hom_{S(n)}
\left(M,M' \otimes \Kosz^\vee\right) 
$$
for any two $S(n)$-modules $M$ and $M'$.

\subsection{}

In particular, the power-sum functions $p_\mu$ satisfy 
$$
\Hom_{S(n)}(p_\mu,M') = \tr_{M'} \sigma_\mu
$$
where $\sigma_\mu$ is a permutation with cycle type $\mu$. The 
trace of $\sigma_\mu$ in $\Kosz$ is easy to compute and we get  
\begin{equation}
  \lang p_\lambda,p_\mu \rang = \delta_{\lambda,\mu} \, \zz(\mu) \, \ww(\mu)
\label{HermPro}
\end{equation}
where $\zz(\mu)$ is the order of the centralizer of the conjugacy 
class $\mu$ and 
$$
\ww(\mu)= \prod_i (q^{\mu_i/2}-q^{-\mu_i/2})
(t^{\mu_i/2}-t^{-\mu_i/2})
$$

\subsection{}

The isomorphism ${\mathcal{B}}$ on the left in \eqref{comm1} is the map constructed
by Haiman and Bridgeland-King-Reid. Without going in the 
details of this construction, it suffices to say that it sends the
skyscraper sheaves at the torus fixed points to Macdonald 
polynomials. 

The torus fixed points of $\Hilb_n$ 
are the ideals of the form
$$
I_\lambda = (x^\lambda_1,x^{\lambda_2}y,x^{\lambda_3}y^2,\dots)
$$
where $\lambda$ is a partition of $n$. The symmetric function corresponding to 
$\cO_{I_\lambda}$ is 
\begin{equation}
  H_\lambda = t^{n(\lambda)} \, 
\Upsilon \, J_\lambda(q,t^{-1}) \,,
\label{defH}
\end{equation}
where $J_\lambda$ is the integral form of the Macdonald polynomial, as defined 
in Macdonald's book \cite{Mac}, $\Upsilon$ is an algebra automorphism 
defined by (cf. \eqref{eq:alk}):
\begin{equation}
  \Upsilon \, p_k = (1-t^{-k})^{-1} p_k \,,\label{defUps}
\end{equation}
and 
$$
n(\lambda) = \sum (i-1) \lambda_i \,.
$$ 

\subsection{}

The map $\mathcal{B}$ comes from an equivalence of derived categories and, 
hence, it is an isometry of $K$-groups. Concretely this means 
that $H_\mu$ are orthogonal with respect to \eqref{HermPro} 
\begin{equation}
  \label{Hm2}
  \lang H_\lambda, H_\mu \rang = \delta_{\lambda,\mu}  \,
\prod_{\textup{cotangent weights $w$}} (w^{1/2}-w^{-1/2})\,,
\end{equation}
where $w$ ranges over all torus weights in the cotangent space to $\Hilb_n$ at $I_\mu$. The $T$-character
of this space is well-known to be 
$$
\sum_{\square\in\lambda} q^{a(\square)+1} t^{-l(\square)}+ q^{-a(\square)} t^{l(\square)+1}\,,
$$
where $a(\square)$ and $l(\square)$ denote the arm and the leg of a square, as usual. 

\subsection{}
The line bundle 
$$
\cO(1) = \Lambda^{\textup{top}} \, \left(\cO/I \right)
$$
is the ample generator of the Picard group of the Hilbert scheme. 
The operator 
$$
\cF \xrightarrow{\,\bT\,} \cF(1) = \cF \otimes \cO(1) 
$$
plays an important role in the theory. By construction 
$$
\bT \, H_\lambda = q^{n(\lambda')} \, t^{n(\lambda)} \, H_\lambda \,.
$$
This eigenvalue is the determinant of the torus action in 
$\cO/I_\lambda$.

\subsection{}

By definition, 
$$
\bV \cO_{I_\mu} = \bigoplus_{k>0} \left(I_\mu\right)^{\boxtimes k} 
$$
Therefore
\begin{equation}
\lang [\lambda] \otimes \cO, \bV H_\mu \rang = (qt)^{|\lambda|/2} \, 
S^\lambda(I_\mu)  
\label{SV}
\end{equation}%
where 
$S^\lambda$ denotes the Schur functor 
$$
S^\lambda(A) = \Hom_{S(n)}
([\lambda],A^{\otimes n}) \,, \quad  n = |\lambda| \,.  
$$

\subsection{}
In particular, the map 
$$
\bL \xrightarrow{\quad \lang \,\cdot\, , \bV H_\mu \rang \quad} 
\bR
$$
is an antilinear algebra homomorphism. That is, 
there is a point $\xmu$ in 
the spectrum of $\bL$ such that 
\begin{equation}
\lang f , \bV H_\mu \rang  = (qt)^{-\deg f/2} \, \bar{f}(\xmu) \,.
\label{pairbV}
\end{equation}
The prefactor here is just to avoid fractional powers later.

\subsection{} 

The point $\xmu$ may be described very concretely. 
{} From definitions
\begin{align}
  p_k(\xmu) &= (qt)^{k/2} \lang p_k, \bV H_\mu \rang \notag \\
&= (1-q^k)(1-t^k)
  \, p_k(I_\mu) \notag  \\
&= 
p_k\left( \left\{\textup{generators of
        $I_\mu$}\right\}\right)
  -p_k\left(\left\{\textup{relations}\right\}\right)\,.
\label{pofmu} 
\end{align}
where 
$p_k(I_\mu)$ is the sum of $k$-th powers of 
of the torus weights in $I_\mu$ and the arguments in 
$p_k\left( \left\{\textup{generators of
        $I_\mu$}\right\}\right)$ are the torus 
weights of generators in a minimal free resolution of $I_\mu$. 

For example, for $I_\square$, we have two generators 
of weights $q,t$ and one relation of weight $qt$, so 
$$
p_k(\xmbox) =  
q^k+t^k- (qt)^k \,. 
$$



\subsection{} 
Now we can state the key result which gives a formula for 
the operator $\bV$ in terms of a Heisenberg algebra 
action on $\bL$. 

We define the creation operators as multiplication by 
$$
\bal_{-n} = \frac{p_n}{(1-q^n)(1-t^n)}\,, \quad n>0\,.
$$
We define the annihilation operators as 
adjoint operators 
$$
\bal_n = \bal_{-n}^* \,.
$$ 
It follows that
$$
[\bal_n,\bal_m] = \delta_{n+m} \, n \, \ww(n)^{-1} \,. 
$$
or, in other words,  
$$
\bal_n = n \, (qt)^{n/2} \, \frac{\partial}{\partial p_n} \,, \quad n>0 \,. 
$$

\subsection{}

We also define the corresponding holomorphic/singular 
parts of a deformed boson 
$$
\bphi_\pm(z) = \sum_{\mp n>0} \frac{\bal_{n}}{n} \, z^{-n} \,.
$$
These satisfy
$$
\bphi_{\pm}(z)^* = - \bphi_{\mp}(1/z) 
$$
and 
\begin{equation}
  \label{comm_phi}
  \left[\bphi_-(z),\bphi_+(w)\right] = \log\,  \left( \sqrt{qt} \, w/z\right)_{q,t} 
\end{equation}
where the function $(x)_{q,t}$ was defined in \eqref{eq:xqt}.

\subsection{}

\begin{Theorem}\label{T2}
$$
\bV = 
(-1)^{L_0} \bT \,  e^{\bphi_+(1)} e^{\bphi_-(\sqrt{qt})}
$$
\end{Theorem}

\noindent 
Note that in Theorem \ref{T2} 
 the operator $\bT$ appears after the vertex operator, that is, it is applied on the orbifold side,
and not on the Hilbert scheme side where it acts naturally. 

Proof of Theorem \ref{T2} will be given in Section \ref{s_proof}. 
 
\subsection{}\label{s_C1}

\begin{Corollary}\label{C1} 
We have 
$$
\bW(m) = \left(\frac{m}{\sqrt{qt}}\right)^{L_0}\!\!
e^{\bphi_+(1)-\bphi_+(qt/m)} \, e^{\bphi_-(\sqrt{qt}) - \bphi_-(\sqrt{qt}/m)} 
$$  
\end{Corollary}

\noindent 
This follows from Theorems \ref{T1} and \ref{T2} and commutation 
relations.



\section{Proof of Theorem \ref{T2}}
\label{s_proof}

\subsection{}

Our strategy is to reduce Theorem \ref{T2} to Macdonald-Mehta-Cherednik identities
for $GL(N)$. In particular, the vertex operators will arise as a stabilization 
of theta functions as $N$ grows to infinity. Here it will be important 
that theta functions for the weight lattice of $GL(N)$ have a factorization, 
namely the Jacobi triple product. 

The MMC identities are best understood in the context of $SL(2,\Z)$ action by 
automorphisms of Cherednik algebras, see 
\cite{C1,C2,Kir}. The same $SL(2,\Z)$ action plays a very 
important role in supersymmetric gauge theories. In the present note, we skip
a proper discussion of this topic and go straight for the formulas.

\subsection{}

We can use the operator \eqref{defUps} 
$$
\Upsilon: \bL \to \bL 
$$
to pull back the Hermitian form 
$$
\lang f,g \rang_\Upsilon = \lang \Upsilon f, \Upsilon g\rang \,. 
$$
Then 
\begin{equation}
  \lang f, \bV^\Upsilon J_\lambda(q,t^{-1}) \rang_\Upsilon = 
t^{-n(\lambda)} \bar{f} (q^{\lambda_1},q^{\lambda_2} t, q^{\lambda_3} t^2, \dots) \label{bBY}
\end{equation}
where 
$$
\bV^\Upsilon = \Upsilon^{-1} \, \bV \, \Upsilon \,.
$$
Indeed, it is enough to check \eqref{bBY} on the generators $p_k$ in which 
case it reduces to \eqref{pofmu}. 

\subsection{}
After replacing $t$ by $t^{-1}$, the Hermitian form 
$\lang \,\cdot\,,\,\cdot\, \rang_\Upsilon$ is the limit of 
the normalized Macdonald inner product for $GL(N)$ as $N\to\infty$, that is
$$
\lang f,g \rang_\Upsilon = \lim_{N\to \infty} \frac{\lang \bar f \, g \rang_N}{\lang 1 \rang_N}
$$
where, cf. \eqref{eq:macdom}, 
$$
\lang g \rang_N  = \frac{1}{N!} \int_{|x_i|=1} g(x) \, \Delta(x) \, \textup{dHaar} 
$$
and, cf. \eqref{eq:macdo1},
$$
\Delta(x) = \prod_{\alpha} \frac{(x^\alpha;q)_\infty}{(x^\alpha/t;q)_\infty} \,.
$$
Here $\alpha$ ranges over all roots of $GL(N)$. 
\subsection{}

In Cherednik theory, it is well-known that the operator that takes 
Macdonald polynomials $J_\lambda$ to evaluation at $q^{\lambda+\rho}$ as 
in \eqref{bBY} it the theta function of the weight lattice. So 
Theorem \ref{T2} becomes, basically, a matter of 
taking the limit of the theta function as $N\to\infty$. For
$GL(N)$, it is easy because of the factorization of the theta function.

\subsection{} 

Cherednik's formula for the root system $A_{N}$ says, in coordinates,
\begin{multline}
  \frac{1}{N!}\int  dx \cdot P_\mu(x) P_\nu(x^{-1})\Theta(x) \Delta(x;q,t) =  \\ 
  q^{\frac{\mu^2}{2}+\frac{\nu^2}{2}}t^{-n(\nu)}
  P_\nu(q^{-\mu-\rho})\prod_{i<j} 
  \frac{(q^{\mu_i-\mu_j}t^{j-i};q)_\infty}
  {(q^{\mu_i-\mu_j }t^{j-i+1};q)_\infty}.
\end{multline}
Here $t=q^k$ for $k$ a positive integer, $N > |\mu|,|\nu|$, and
\[f(q^{-\mu-\rho}) = f(q^{-\mu_1},q^{-\mu_2}t,q^{-\mu_3}t^2,...), \quad \mu^2 = \sum_i \mu_i^2.\]
Our notation differs from formula 1 of \cite{C2}, in that
we use the parameters $q,t$, for the Macdonald polynomials instead of $q^2,t^2$,
and
\[\Delta(x;q,t) = \prod_{i\neq j} \frac{(x_ix_j^{-1};q)_\infty}{(x_ix_j^{-1}t;q)_\infty}.\]
Secondly, our definition of $f(q^\nu)$ differs by a factor of $t^{-N/2}$, which is canceled by
$t^{\langle \rho,\nu\rangle-n(\nu)}$, and $\Theta$ is the theta function on the $\Z^N$ lattice,
$$
\Theta = \prod_{i=1}^N \vartheta(x_i;q)\,, \quad \vartheta(x) = \sum_n x^n q^{n^2/2} \,.
$$
Using the full lattice does not affect the integral since $N$ is large enough.

By the Jacobi triple product formula,
$$
\Theta(x) = (q;q)_\infty^N \, \Theta_+(x) \, \Theta_+(x^{-1}), 
$$
where
$$
\Theta_+ = \exp\left(\sum_{k>0} \frac{(-1)^k}{q^{k/2}-q^{-k/2}} \, p_k \right).
$$
Now divide the left hand side by $c_N(q;q)_\infty^N$, where $c_N$ is given by  (\cite{Mac} chapter VI, (9.7)),
\[ c_N = \frac{1}{N!} [\Delta]_1\]
and $[\Delta]_1$ denotes the constant term in $x_i$.
Taking the limit
as $N\to\infty$ gives
\[
  \left(\omega \Gamma_-\left(\frac{q^{1/2}}{1-q}\right)\omega \cdot P_\mu,
  \omega \Gamma_-\left(\frac{q^{1/2}}{1-q}\right)\omega \cdot P_\nu\right)_{q,t},\]
where 
\begin{equation}
\omega \, p_k = (-1)^{k-1} p_k 
\end{equation}
is the standard involution on symmetric functions and 
\[\left(p_\mu,p_\nu\right)_{q,t} = \delta_{\mu,\nu}\zz(\mu)\prod_k \frac{1-q^{\mu_k}}{1-t^{\mu_k}}.\]
is the standard Macdonald symmetric inner-product. 

Using the commutation relations, we have
\begin{multline}
\label{finmac}
  \left(\omega \Gamma_-\left(\frac{q^{1/2}}{1-q}\right)\Gamma_+\left(\frac{q^{1/2}}{1-t}\right)
  \omega \cdot P_\mu,P_\nu\right)_{q,t} =  \\ 
  q^{\frac{\mu^2}{2}+\frac{\nu^2}{2}}t^{-n(\nu)}
  P_\nu(q^{-\mu-\rho})\prod_{i} \left(\frac{(t^i;q)_\infty}{(t;q)_\infty}\prod_{j>i}
  \frac{(q^{\mu_i-\mu_j}t^{j-i};q)_\infty}
  {(q^{\mu_i-\mu_j }t^{j-i+1};q)_\infty} \right)= \\
  q^{\frac{\mu^2}{2}+\frac{\nu^2}{2}}t^{-n(\nu)}
  P_\nu(q^{-\mu-\rho}).
\end{multline}

For the remainder of the paper define a symmetric inner-product by
\[(p_\mu,p_\nu)_{q,t}' = \delta_{\mu,\nu}\mathfrak{z}(\mu)\prod_{k} (1-q^{\mu_k})(1-t^{\mu_k}).\]
Given an expression $f$, let $f_{[k]}$ denote the result of substituting $x^k$ in for $x$ for all arguments in $f$,
and define the vertex operators on $\Lambda$ by
\[\Gamma_-(f) = \exp\left(\sum_{k>0} \frac{f_{[k]}}{k} p_k\right),\quad
\Gamma_+(f) = \exp\left(\sum_{k>0} f_{[k]}\frac{\partial}{\partial p_k}\right),\]
so that 
\[e^{\bphi_+(z)} = \Gamma_-\left(-\frac{z}{(1-q)(1-t)}\right),\quad 
e^{\bphi_-(z)} = \Gamma_+\left(z^{-1}\sqrt{qt}\right)\]

Now we find the matrix elements of the right hand side in theorem \ref{T2},
\[\lang H_\mu, (-1)^{L_0} \bT \,  e^{\bphi_+(1)} e^{\bphi_-(\sqrt{qt})}H_\nu \rang =\]
\[\pm \left( \omega\overline{ T^{-1} H_\mu}, (qt)^{-L_0/2} e^{\bphi_+(1)} e^{\bphi_-(\sqrt{qt})}H_\nu \right)'_{q,t} =\]
\[\pm q^* t^*  \left( H_\mu,  e^{\bphi_+(1)} e^{\bphi_-(\sqrt{qt})} H_\nu \right)_{q,t}' =\]
\[\pm q^* t^* \left( \Upsilon J_\mu(q,t^{-1}) , 
 \Gamma_-\left(-\frac{1}{(1-q)(1-t)}\right) \Gamma_+(1)\Upsilon   J_\nu(q,t^{-1}) \right)'_{q,t} =\]
\[\pm q^* t^* \left( \Upsilon J_\mu(q^{-1},t) , \Upsilon
\Gamma_-\left(\frac{t^{-1}}{1-q}\right)\Gamma_+\left(\frac{1}{1-t^{-1}}\right) J_\nu(q^{-1},t) \right)'_{q,t} =\]
\[\pm  q^*t^* \left.\left( J_\mu ,
\Gamma_-\left(\frac{t^{-1}}{1-q^{-1}}\right)\Gamma_+\left(\frac{1}{1-t^{-1}}\right) J_\nu \right)_{q,t}\right|_{q=q^{-1}} =\]
\begin{equation}
\label{thm2eq}
\pm  q^*t^* \left.\left( J_\mu ,
\omega \Gamma_-\left(\frac{q^{1/2}}{1-q}\right)\Gamma_+\left(\frac{q^{1/2}}{1-t}\right) \omega J_\nu \right)_{q,t}\right|_{q=q^{-1}}.
\end{equation}
Here we have left out the sign and the powers of $q,t$ outside the expression for simplicity, and have made use of the following facts about symmetric polynomials,
\[\Gamma_{\pm}(f) x^{L_0} = x^{L_0} \Gamma_{\pm}(fx^{\pm1}),\quad 
\Gamma_{\pm}(f) \omega = \omega (-1)^{L_0}\Gamma_{\pm}(-f) (-1)^{L_0}, \]
\[(\Upsilon f,\Upsilon g) =  t^{2\text{deg}(f)} (f,g)_{q,t},\quad  J_\mu(q,t^{-1}) =  q^*t^* J_\mu(q^{-1},t).\]
The last formula follows from the fact that $P_\mu$ is real, $\overline{P_\mu} = P_\mu$.

Inserting the integral form
\[J_\mu = \prod_{\Box \in \mu} (1-q^{a(\Box)}t^{\ell(\Box)+1}) P_\mu,\]
into Eq. \eqref{thm2eq}, and using the formula \eqref{finmac}, we get
 \[\pm q^*t^* \left. J_\mu(q^{-\nu-\rho}) \right|_{q=q^{-1}}.\]
By the definition of $H_\mu$, we have
\[\overline{H_\mu}(x_\nu) =\left. t^* J_\mu(q^{-\nu-\rho})\right|_{q=q^{-1}}.\]
Using this relation, keeping track of the signs and powers of $q,t$, we end up with
\[\lang H_\mu, \bV H_\nu\rang,\]
from which the theorem follows.
\qed

\section{Interpolation polynomials}\label{s_interp}

\subsection{} 

In this section, we interpret our results in terms of 
interpolation Macdonald polynomials,  studied by by the G.~Olshanski and 
one of the authors, F.~Knop, S.~Sahi, E.~Rains, and many others, see for example
 \cite{KS, OO, Rains} and also \cite{Four} for an elementary introduction to the subject. 
Interpolation Macdonald polynomials are remarkable inhomogeneous symmetric functions defined by Newton-style interpolation conditions. 
 As it turns out, 
the original orthogonal Macdonald polynomials are their terms of top degree. 

Macdonald polynomials generalize characters of semisimple Lie groups, that is, 
of $G=GL(N)$ in the case at hand, 
while Macdonald-Ruijsenaars commuting difference operators 
\cite{Rui} generalize the center 
$\mathcal{Z}\left(\mathcal{U}\mathfrak{g}\right)$ of the 
universal enveloping algebra of $G$, or more precisely 
the radial parts of the invariant differential operators on $G$. 
There is a natural Fourier pairing 
\begin{equation}
\mathcal{Z}\left(\mathcal{U}\mathfrak{g}\right) \otimes 
\C[G]^{G} \to \C \label{FourierG}
\end{equation}
between differential operators and functions on the 
group that applies an operator to a function and evaluates the result at a 
special point, in this case, the identity element. 

This may be used to construct 
a distinguished linear basis of 
$\mathcal{Z}\left(\mathcal{U}\mathfrak{gl}(N)\right)$ formed by elements
of a given order that vanish in as many irreducible representations as possible, 
see \cite{OO1} for a comprehensive discussion. The Harish-Chandra isomorphism 
sends this basis to the basis of interpolation Schur functions $s^*_\lambda$, 
that becomes the basis $P^*_\mu$ of interpolation Macdonald polynomials under
the deformation. 

\subsection{}

Concretely, the interpolation Macdonald polynomials $P^*_\mu$ are the unique, 
up to multiple, symmetric polynomials of degree $|\mu|$ such that 
\begin{equation}
P^*_\mu(q^\lambda t^\rho) = 0 \,, \quad  \mu \not\subset \lambda\,,
\label{Pvan}
\end{equation}
where $q^\lambda t^\rho = (q^{\lambda_1},q^{\lambda_2} t^{-1},q^{\lambda_3} t^{-2},
\dots)$. These can be normalized to be monic, and, in fact, to satisfy
$$
P^*_\mu = P_\mu + \dots 
$$
where $P_\mu$ are the monic orthogonal Macdonald polynomials and dots stand for 
terms of lower degree. 

After the transformations that take $P_\mu$ to $H_\mu$, the evaluation condition
in \eqref{Pvan} becomes 
$$
H^*_\mu(\ominus \mathsf{x}_\lambda) = 0 \,, 
$$
where $p_k(\ominus \mathsf{x}_\lambda) = - p_k(\mathsf{x}_\lambda)$. In formulas, 
this minus sign appears as the mismatch between the denominator in 
$$
p_k(q^{\lambda} t^{-\rho}) = \frac{1}{1-t^k} \, 
\left(
p_k\left( \left\{\textup{generators of
        $I_\lambda$}\right\}\right)
  -p_k\left(\left\{\textup{relations}\right\}\right)\right)
$$
and the similar, but not identical factor in \eqref{defUps}. 

\subsection{}
The pairing \eqref{FourierG} becomes the following Fourier pairing 
\begin{equation}
\left(f,H_\mu\right)_\textup{Fourier} = f(\ominus \xmu) \, H_\mu(\ominus 
\mathsf{x}_\varnothing) \,.
\label{FourM}
\end{equation}
It is easy to relate it to the Hermitian pairing \eqref{pairbV}. A basic
property of interpolation Macdonald polynomials is their orthogonality 
$$
\left(H^*_\mu,H^*_\lambda\right)_\textup{Fourier} = 0 \,, \quad \mu\ne \lambda\,,
$$
with respect to \eqref{FourM}, see \cite{Four}. Formally, this orthogonality 
is equivalent to the binomial formula of \cite{Obin}. The orthogonality 
implies, in particular, the symmetry of Fourier pairing \eqref{FourM}, which is one 
of the cornerstones of Macdonald-Cherednik theory, see e.g.\ \cite{C1}.

The essence of Macdonald-Mehta-Cherednik identities, and also of Theorem 
\ref{T2}, is that the multiplication by the theta function of the 
weight lattice takes Macdonald Hermitian inner product to a Fourier-type 
pairing. One can then use a Wiener-Hopf factorization of the theta function, 
or equivalently 
the factorization in Theorem \ref{T2} into raising and lowering operators,
to produce an operator that takes an orthogonal basis for one form to  
an orthogonal basis for another. This means it takes 
orthogonal Macdonald polynomials to 
interpolation Macdonald polynomials. Concretely, we have the following

\begin{Proposition}\label{NGa}
\begin{equation}
  \label{GammaN}
  H^*_\mu = \bT^{-1} \, e^{\sum_{n>0} (-1)^{n+1} \frac{\partial}{\partial p_n}} \, \bT \, H_\mu \,. 
\end{equation}
\end{Proposition}

Conversely, Theorem \ref{T2} may be deduced from \eqref{GammaN}. In fact, this was 
our original proof, before we realized the equivalence to the MMC identity. 
Other statements in the literature can be seen equivalent to \eqref{GammaN}, 
in particular Proposition 3.5.17 in \cite{Hai}, which goes back to \cite{GHT}. 

\subsection{}

From this point of view, the form of the vanishing condition in \eqref{Pvan} 
\begin{equation*}
P^*_\mu(q^\lambda t^\rho) = 0 \,, \quad \lambda \notin \subset \mu + 
(\Z_{\ge 0})^\infty \,,
\label{Pvan2}
\end{equation*}
may be traced to the support 
\begin{equation}
(\Z_{\ge 0})^N \subset \Z^N = \textup{weight lattice of $GL(N)$} 
\label{suppTh}
\end{equation}
of the Wiener-Hopf, that is, triple-product factorization of the theta-function. 

Here, by a Wiener-Hopf factorization of a multivariate function we mean 
its factorization as a function of one of its arguments. Thus, for example, 
the support of a Wiener-Hopf factorization of a theta-function $\Theta_L$ of 
a general lattice $L$ is a just some half-space of $L$. A much smaller support 
in \eqref{suppTh} is a special feature of $GL(N)$, and also of 
$BC_n$, which makes
the theory of interpolation polynomials much richer for these root system. 

This explains why one has a very good theory of interpolation polynomials
for $GL(N)$ but not, for example, for $SL(N)$. 

\section{Applications}

One application of our result is the calculation of the partition $\cZ$-functions
of the five dimensional $A_r$-type quiver $U(1)^{r+1}$ gauge theories. 
The latter are characterized by $r+1$ coupling constants
${\qe}_{0}, {\qe}_{1}, \ldots , {\qe}_{r}$, and $r+1$ masses ${\mu}_{0}, {\mu}_{1}, {\mu}_{2}, \ldots , {\mu}_{r}$. 

In this section we change the notations $(q,t) \mapsto (q_{1}, q_{2})$. 
The partition function is given by:
\begin{equation}
Z( {\mu}, {\qe}; q_{1}, q_{2}) = {\rm Tr}\, \prod_{i=0}^{r} \, {\qe}_{i}^{L_{0}} \, {\bW}(m_{i}) \,
\label{eq:zfquiv}
\end{equation}
We should stress that \eqref{eq:zfquiv} is obtained using the five dimensional gauge theory
considerations. The formalism described in this paper allowed to repackage the instanton sum into the trace of a product of vertex operators, a typical expression for the correlation function of the chiral operators in a $q$-deformed conformal field theory on a torus.  
The expression \eqref{eq:zfquiv} can be now used to test the duality between the IIA string theory and the M-theory, in that the gauge theory we are discussing in five dimensions
is obtained by twisted compactification of the $(2,0)$-theory with some defects on an elliptic curve. Let us discuss this in some more detail now.

Note that for $r=0$ the theory we are discussing is the $U(1)$ ${\cN}=2^{*}$ theory, i.e. the theory with the adjoint hypermultiplet. The latter can be viewed geometrically as the compactification of the $(2,0)$ theory on the elliptic curve with twisted boundary conditions on both the worldvolume and the transverse space of the theory (i.e. the twist involves the $R$-symmetry group $Spin(5)$). The noncompact part of the worldvolume, the Euclidean space ${\R}^{4} \approx {\C}^{2}$, with the coordinates $1234$ is twisted with the parameters $(q_{1}, q_{2})$, while the transverse ${\R}^5 \approx {\C}^{2} \oplus {\R}^1$, 
with the coordinates $56789$ is twisted with the parameters $( m / \sqrt{q_{1}q_{2}}, m^{-1} / \sqrt{q_{1}q_{2}} ) \oplus 1$. 
The trace calculation \eqref{eq:zfquiv} for $r=0$ gives, cf. \cite{NO}:
\begin{equation}
Z^{\rm inst}({\mu}, {\qe}; q_{1}, q_{2}) = {\rm Tr}\, {\qe}^{L_{0}} {\bW}(m) = 
{\exp}\, \sum_{n=1}^{\infty} \frac{Q^{n}}{n\, m^n} \frac{(m^n - q_{1}^n)(m^n - q_{2}^{n})}{(1-Q^n) ( 1 - q_{1}^{n}) ( 1- q_{2}^{n}) } 
\label{eq:zfunrz}
\end{equation}
where $$
Q = \frac{{\qe} m }{\sqrt{q_{1}q_{2}}}
$$
In this way we prove the conjecture of \cite{NJ} modulo the redefinition $m \mapsto m q_{1}q_{2}$. Actually, the expression \eqref{eq:zfunrz}
does not compute the full partition function of the five dimensional gauge theory. Indeed, the latter must agree with the partition function of the $(2,0)$ theory compactified on a torus, and enjoy some modular properties, as a consequence. It is clear what is the missing piece. The formula \eqref{eq:zfunrz} is obtained by analyzing the effects of instantons, yet there is also a purely perturbative contribution, which is $\qe$-independent. 

It is not hard to compute it, with the result:
\begin{equation}
Z^{\rm pert} ({\mu}; q_{1}, q_{2}) =  {\exp}\, \sum_{n=1}^{\infty} \frac{q_{1}^{n}q_{2}^{n}}{n} \frac{1-m^n}{ ( 1 - q_{1}^{n}) ( 1- q_{2}^{n}) }
\label{eq:zfunpert}
\end{equation}
Then $Z = Z^{\text{pert}}Z^{\text{inst}}$ is equal to the partition function 
\begin{equation}
Z = {\text{Tr}}_{\cH} \left( (-1)^{\bar F} {\qe}^{L_{0}} {\bar\qe}^{{\bar L}_{0}} \ q_{1}^{J_{12}-\frac 12 (J_{56}+J_{78}) } q_{2}^{J_{34}-\frac 12 (J_{56}+J_{78})} m^{J_{56} - J_{78}} \right)
\label{eq:trh}
\end{equation}
of the free $(2,0)$ six dimensional tensor multiplet computed as a trace over the Hilbert space $\cH$ obtained by quantization on ${\bf S}^{1} \times {\R}^{4}$. It would be interesting to extend this analysis to the case of general $A_{r}$ quiver theories.

\end{document}